\documentclass{article}
\usepackage[T1]{fontenc}
\usepackage{amsthm, fullpage}
\usepackage{amsmath}
\usepackage{amssymb}
\usepackage{amsfonts}
\usepackage{eucal}
\usepackage[all]{xy}
\usepackage[frenchb, english]{babel}
\usepackage{pslatex}
\usepackage[pagebackref]{hyperref}

\newtheorem{prop}{Proposition}[subsection]
\newtheorem{theo}[prop]{Théor\`eme}
\newtheorem*{theo**}{Théorème}
\newtheorem{coro}[prop]{Corollaire}
\newtheorem*{conj*}{Conjecture}
\newtheorem{lemm}[prop]{Lemme}
\newtheorem{lemm*}{Lemme}[prop]

\theoremstyle{definition}

\newtheorem{vide}[prop]{}
\newtheorem{defi}[prop]{Définition}
\newtheorem*{defi*}{Définition}

\theoremstyle{remark}
\newtheorem{rema}[prop]{Remarques}
\newtheorem{exem}[prop]{Exemples}

\numberwithin{equation}{prop}

\newcommand{\riso}{ \overset{\sim}{\longrightarrow}\, }

\newcommand{\Spf}{\mathrm{Spf}\,}

\renewcommand{\sp}{\mathrm{sp}}

\renewcommand{\AA}{{\mathcal{A}}}

\newcommand{\E}{{\mathcal{E}}}

\renewcommand{\H}{{\mathcal{H}}}

\newcommand{\D}{{\mathcal{D}}}

\newcommand{\PP}{{\mathcal{P}}}

\renewcommand{\O}{{\mathcal{O}}}

\newcommand{\V}{\mathcal{V}}

\newcommand{\Y}{\mathcal{Y}}

\newcommand{\X}{\mathfrak{X}}

\newcommand{\U}{\mathfrak{U}}

\newcommand{\DD}{\mathbb{D}}
\renewcommand{\L}{\mathbb{L}}
\newcommand{\R}{\mathbb{R}}
\newcommand{\Q}{\mathbb{Q}}
\newcommand{\Z}{\mathbb{Z}}

\newcommand{\hdag}{  \phantom{}{^{\dag} }    }

\begin{document}
\title{Une caractérisation de la surcohérence}
\author{Daniel Caro \footnote{L'auteur a bénéficié du soutien du réseau européen TMR \textit{Arithmetic Algebraic Geometry}
(contrat numéro UE MRTN-CT-2003-504917).}}
\date{}

\maketitle

\begin{abstract}
\selectlanguage{english}
  Let $\mathcal{P}$ be a proper smooth formal $\mathcal{V}$-scheme, $X$ a closed subscheme of the special fiber of $\mathcal{P}$,
$\mathcal{E} \in F\text{-}D ^\mathrm{b} _\mathrm{coh} ( \D ^\dag _{\mathcal{P},\mathbb{Q}})$ with support in $X$.
We check that $\mathcal{E}$ is $\D ^\dag _{\mathcal{P},\mathbb{Q}}$-overcoherent if and only if,
for any morphism $f\,:\, \mathcal{P}' \rightarrow \mathcal{P}$ of smooth formal $\mathcal{V}$-schemes,
$f ^! (\mathcal{E}) $ is $\D ^\dag _{\mathcal{P}', \, \mathbb{Q}}$-coherent.
\end{abstract}

\selectlanguage{frenchb}

\tableofcontents

%\selectlanguage{frenchb}
%\begin{abstract}
%   Soient $\mathcal{P}$ un $\mathcal{V}$-schéma formel propre et lisse, $X$ un sous-schéma fermé de
%la fibre spéciale de $\mathcal{P}$,
%$\mathcal{E} \in F\text{-}D ^\mathrm{b} _\mathrm{coh} ( \D ^\dag _{\mathcal{P},\mathbb{Q}})$ à support dans $X$.
%On vérifie que $\mathcal{E}$ est $\D ^\dag _{\mathcal{P},\mathbb{Q}}$-surcohérent si et seulement si,
%pour tout morphisme $f\,:\, \mathcal{P}' \rightarrow \mathcal{P}$ de $\mathcal{V}$-schémas formels lisses,
%$f ^! (\mathcal{E}) $ est à cohomologie $\D ^\dag _{\mathcal{P}', \, \mathbb{Q}}$-cohérente.
%\end{abstract}

\section*{Introduction}
Nous aimerions disposer d'une \textit{bonne} cohomologie $p$-adique, i.e., d'une catégorie de $\Q _p$-objets sur des variétés de caractéristique $p$
stable par les six opérations cohomologiques de Grothendieck, à savoir $\otimes$,  $\mathcal{H}om$, $f _*$, $f ^*$, $f _!$ et $f^!$.
Pour cela, en s'inspirant de la caractéristique nulle,
Berthelot a élaboré une théorie des \textit{$\D$-modules arithmétiques} (pour une introduction, on consultera \cite{Beintro2}).
Plus précisément,
soit $\mathcal{V}$ un anneau de valuation discrète complet d'inégales caractéristiques $(0,p)$,
de corps des fractions $K$, de corps résiduel $k$ supposé parfait.
Soit $\PP$ un $\V$-schéma formel lisse.
Il construit (voir \cite{Beintro2} ou \cite{Be1}) alors le faisceau des opérateurs différentiels de niveau fini sur $\PP$,
noté $\D ^\dag _{\PP, \, \Q}$.
Ce faisceau est le complété faible $p$-adique (d'où le symbole {\og $\dag$\fg})
tensorisé par $\Q$ (d'où le symbole {\og $\Q$\fg})
du faisceau $\D _{\PP}$ des opérateurs différentiels usuel sur $\PP$.
Un $\D$-module arithmétique sur $\PP$ signifie un
$\D ^\dag _{\PP, \, \Q}$-module (toujours à gauche par défaut).
Un $F\text{-}\D ^\dag _{\PP, \, \Q}$-module (resp. $F$-objet) est
un $\D ^\dag _{\PP, \, \Q}$-module (resp. objet) muni d'une structure de Frobenius.
En outre, il a défini l'holonomie de manière analogue au cas classique :
un $F\text{-}\D ^\dag _{\PP, \, \Q}$-module cohérent est holonome s'il est nul ou si la dimension de sa variété
caractéristique est égale à la dimension de $\PP$.
Rappelons que pour définir canoniquement la variété caractéristique, on a besoin de la structure de Frobenius
afin de descendre au niveau $0$ (surtout, sans la structure de Frobenius, la variété caractéristique dépend a priori du niveau choisi :
voir \cite[5.2.4]{Be2}).
L'holonomie est une notion stable par foncteur dual (voir \cite[III.4.4]{virrion}) et par produit tensoriel externe (voir \cite[5.3.5.(v)]{Beintro2}).
Concernant les autres opérations, il conjecture (voir \cite[5.3.6.]{Beintro2}) que la catégorie des $F$-complexes
de $\D$-modules arithmétiques à cohomologie bornée et holonome est stable par image directe par un morphisme propre,
par image inverse extraordinaire (ces deux conjectures impliquant d'ailleurs la stabilité de l'holonomie
par image directe extraordinaire, image inverse et par foncteur cohomologique local).

Afin de disposer d'une catégorie stable par de telles opérations
cohomologiques, nous avions dans \cite[3]{caro_surcoherent} défini la notion de $\D ^\dag _{\PP, \, \Q}$-surcohérence.
Par définition, un complexe de $\D ^\dag _{\PP, \, \Q}$-modules
à cohomologie bornée et cohérente est surcohérent si sa cohérence est préservée
par foncteur cohomologique local et image inverse extraordinaire par un morphisme lisse.
Nous avions vérifié que cette notion de surcohérence est préservée par l'image directe d'un morphisme propre, par image
inverse extraordinaire et par foncteur cohomologique local.
Dans ce travail, pour tout $\V$-schéma formel $\PP$ propre et lisse,
nous établissons la caractérisation suivante de la $\D ^\dag _{\PP, \, \Q}$-surcohérence qui n'utilise pas le foncteur cohomologique local :
un $F$-complexe $\E$ de
$\D ^\dag _{\PP, \, \Q}$-modules à cohomologie bornée et cohérente
est $\D ^\dag _{\PP, \, \Q}$-surcohérent si et seulement si,
pour
 tout morphisme
$f\,:\, \PP' \rightarrow \PP$ de $\V$-schémas formels lisses,
$f ^! (\E) $ est à cohomologie $\D ^\dag _{\PP', \, \Q}$-cohérente (voir \ref{carac-surcoh}).

\medskip

Voyons maintenant le contenu des deux parties de cet article.
Dans la première partie, nous définissons la notion de surcohérence sur un sous-schéma fermé $X$ de la fibre spéciale de $\PP$ de la façon suivante :
un complexe de $\D ^\dag _{\PP, \, \Q}$-modules est surcohérent sur $X$ si,
pour tout morphisme $f\,:\, \PP' \rightarrow \PP$ de $\V$-schémas formels lisses tel que $f (\PP') \subset X$,
$f ^! (\E) $ est à cohomologie $\D ^\dag _{\PP', \, \Q}$-cohérente (voir \ref{defi-SC}).
Nous vérifions de plus que pour tout $F$-complexe $\E$ surcohérent sur chacun de ses points fermés,
il existe un ouvert lisse dense de son support sur lequel les espaces de cohomologie de $\E$
soient associés à des $F$-isocristaux convergents.
Les $\D$-modules arithmétiques provenant d'isocristaux (sur)convergents ont
un énorme avantage : il est possible de faire de la descente propre génériquement finie et étale, ce qui permet en l'occurrence
d'utiliser le théorème de désingularisation de de Jong.

Dans la seconde partie, nous établissons que lorsque $\PP$ est propre,
un $F$-complexe de $\D ^\dag _{\PP, \, \Q}$-modules $\E$ est $\D ^\dag _{\PP, \, \Q}$-surcohérent si et seulement s'il est surcohérent sur son support.
L'idée est de procéder par récurrence sur le support de $\E$.
En gros, notamment en utilisant le fait que $\E$ est aussi surcohérent sur ses points fermés,
il devient possible de dévisser $\E$ en un $F$-complexe dont
le support est de dimension strictement inférieure à celle de $\E$
et en un $F$-complexe dont les espaces de cohomologie
sont associés à des $F$-isocristaux surconvergents via l'équivalence de
catégories usuelle (voir \ref{eqcat-gen}).
En procédant par récurrence sur la dimension du support de $\E$,
il suffit alors d'établir le résultat pour les $F$-complexes dont les espaces de cohomologie
sont associés à des $F$-isocristaux surconvergents, ce qui s'établit par descente via le théorème de désingularisation de de Jong.

\begin{center}
{\bf Notations}
\end{center}

Tout au long de ce travail, nous garderons les notations
suivantes :
les schémas formels seront notés par des lettres calligraphiques ou
gothiques et leur fibre spéciale par les lettres romanes
correspondantes. Si $f\,:\, X' \rightarrow X$ est un morphisme de schémas ou de schémas formels, on note
$d _X$ la dimension de $X$ et $d _{X'/X}$ la dimension relative de $f$.
De plus, la lettre $\V$ désignera un anneau de valuation discrète complet,
de corps résiduel parfait $k$ de caractéristique $p>0$, de corps des
fractions $K$ de caractéristique $0$.
On fixe $s\geq 1$ un entier naturel et $F$ sera la puissance
$s$-ième de l'endomorphisme de Frobenius.
Si $\E$ est un faisceau abélien, $\E _\Q$ désignera $\E \otimes _\Z \Q$.
Lorsque cela n'est pas précisé,
les modules sont des modules à gauche.
Les $k$-variétés sont
les $k$-schémas séparés et de type fini.
Tous les $k$-schémas seront {\it réduits}.
Lorsque le symbole ensemble vide {\og$\emptyset$\fg} apparaît dans une notation, nous ne l'indiquerons pas
(e.g., {\og $\D ^\dag  _{\PP,\Q}$\fg} à la place de {\og $\D ^\dag   _{\PP} ( \emptyset ) _\Q$ \fg}).

Si $\AA$ est un faisceau d'anneaux,
$D ^\mathrm{b}( \AA )$
la catégorie des complexes de $\AA$-modules à cohomologie bornée.
De plus, les indices
{\og $\mathrm{coh}$\fg}, {\og $\mathrm{surcoh}$\fg}, {\og $\mathrm{parf}$\fg}
signifient respectivement
{\og  cohérent\fg}, {\og  surcohérent\fg}
{\og  parfait\fg}.

\section{Surcohérence et surcohérence sur un sous-schéma fermé}

\subsection{Rappels et notations}
\label{rapnot}
Pour une introduction efficace aux $\D$-modules arithmétiques de Berthelot, on pourra consulter \cite{Beintro2}.
En ce qui concerne les isocristaux surconvergents de Berthelot, on pourra se reporter à l'ouvrage \cite{LeStum-book-isoc}
pour une étude complète de leur construction.
Ensuite, pour les liens entre ces deux théories, on lira \cite[4]{Be1}, \cite[4.6]{Be2},
\cite{caro-construction}, \cite{caro_devissge_surcoh}, \cite{caro-2006-surcoh-surcv}.
Précisons tout de même nos notations et donnons quelques propriétés récentes, utiles pour la suite.

$\bullet$ Soient $f\,:\, \PP' \rightarrow \PP$ un morphisme de $\V$-schémas formels lisses,
$T$ un diviseur de $P$ tel que $ f ^{-1} (T)$ soit un diviseur de $P'$.
On dispose alors d'un foncteur {\og image directe par $f$\fg} noté
$f_+ \,:\,(F\text{-})D ^\mathrm{b} _{\mathrm{coh}} (\D ^\dag _{\PP'} (\hdag T') _{\Q}) \rightarrow
(F\text{-})D (\D ^\dag _{\PP} (\hdag T) _{\Q})$,
le symbole {\og $(F\text{-})$\fg} signifiant que ce foncteur est valable avec ou sans structure de Frobenius.
On bénéficie de plus du foncteur {\og image inverse extraordinaire par $f$\fg} noté
$f^! \,:\,(F\text{-})D ^\mathrm{b} _{\mathrm{coh}} (\D ^\dag _{\PP} (\hdag T) _{\Q}) \rightarrow
(F\text{-})D (\D ^\dag _{\PP'} (\hdag T') _{\Q})$. On note aussi $f ^* := \H ^0 ( f^! [-d_{P'/P}])$.
Lorsque $f$ est propre, le foncteur $f _+$ préserve la cohérence.
Lorsque $f$ est lisse, le foncteur $f ^!$ préserve la cohérence et $f ^* \riso  f^! [-d_{P'/P}]$.

$\bullet$ Lorsque le morphisme $f$ est propre, la version \cite[1.2.7]{caro_courbe-nouveau} (i.e., {\og $D ^\mathrm{b} _\mathrm{coh}$\fg}
au lieu de {\og $\widetilde{D} _\mathrm{parf}$\fg}
qui est une sous-catégorie pleine de $D _\mathrm{parf}=D ^\mathrm{b} _\mathrm{coh}$ : voir les notations du début de \cite{Vir04})
du théorème de dualité relative
implique le fait suivant :
pour tous $\E' \in  D _{\mathrm{coh}} ^{\mathrm{b}} ( \D ^{\dag} _{\PP '} (\hdag T')_{\Q} )$,
$\E \in D _{\mathrm{coh}} ^{\mathrm{b}} ( \D ^{\dag} _{\PP } (\hdag T)_{\Q} )$,
on dispose alors de l'isomorphisme canonique d'adjonction
fonctoriel en $\E$ et $\E'$ :
\begin{gather}
\notag
\mathrm{Hom} _{ \D ^{\dag} _{\PP  }( \hdag T)  _{\Q} }( f _{ +} (\E'),\E )
\riso
\mathrm{Hom} _{ \D ^{\dag} _{\PP ' }( \hdag T')  _{\Q}  }(\E'  ,f  ^! (\E) ),
\end{gather}
où $\mathrm{Hom} _{ \D ^{\dag} _{\PP  ,\,\Q}} (-,-) :=
H ^0 \circ \R \mathrm{Hom} _{ D( \D ^{\dag} _{\PP  ,\,\Q})} (-,-) =
\mathrm{Hom} _{D( \D ^{\dag} _{\PP  ,\,\Q})} (-,-)$.
En voici une application (qui nous sera plusieurs fois utile) : lorsque $f$ est propre et lisse,
puisque $f  ^! (\E) $ est alors cohérent, on dispose alors du morphisme d'adjonction :
$f _+ f ^! (\E) \rightarrow \E$.

$\bullet$ Nous exposons brièvement ici les résultats de \cite[2.2]{caro_surcoherent}
concernant le foncteur cohomologique local.
Soit $X$ un sous-schéma fermé d'un $\V$-schéma formel lisse $\PP$.
Le {\og foncteur cohomologique local à support strict dans $X$\fg} de \cite[2.2.6]{caro_surcoherent}
sera noté $\R \underline{\Gamma} ^\dag _X$. Le {\og foncteur de localisation en dehors de $Z$\fg}
(ou abusivement {\og foncteur restriction en dehors de $Z$\fg}) noté $(\hdag X)$ s'inscrit par définition dans le triangle distingué
de localisation en $X$ :
\begin{equation}
  \label{tri-loc-gen}
\R \underline{\Gamma} ^\dag _X (\E) \rightarrow \E \rightarrow \E (\hdag X) \rightarrow \R \underline{\Gamma} ^\dag _X (\E) [1],
\end{equation}
où $\E (\hdag X) $ désigne $(\hdag X) (\E ) $ (les deux notations sont valables).
Ces deux foncteurs $\R \underline{\Gamma} ^\dag _X$ et $(\hdag X)$ vérifient toutes les propriétés attendues, e.g.,
formules classiques de composition, commutation aux images inverses extraordinaires et images directes, triangles distingués de Mayer-Vietoris etc
(voir \cite[2.2]{caro_surcoherent}).

\subsection{Définitions et premières propriétés}
On désigne par $\PP$ un $\V$-schéma formel lisse.
On dispose de la notion de ($F$-)$\D ^\dag _{\PP,\Q}$-surcohérence de \cite[3.1.1]{caro_surcoherent} :
\begin{defi}
\label{defi-311surcoh}
Soit $\E$ un ($F$-)$\D ^\dag _{\PP,\Q}$-module cohérent (resp. un
objet de ($F$-)$D ^\mathrm{b} _{\mathrm{coh}} (\D ^\dag _{\PP,\Q})$). On dit que $\E$ est
{\og $\D ^\dag _{\PP,\Q}$-surcohérent\fg} ou {\og surcohérent sur $\PP$\fg} si pour tout morphisme lisse de $\V$-schémas
formels lisses $f$ : $\PP '\rightarrow \PP$, pour tout
diviseur $T '$ de $P'$, $(\hdag T ') ( f ^* \E)$ est un
($F$-)$\D ^\dag _{\PP',\Q}$-module cohérent (resp. un
objet de ($F$-)$D ^\mathrm{b} _{\mathrm{coh}} (\D ^\dag _{\PP',\Q})$).
On note ($F$-)$D ^\mathrm{b} _{\mathrm{surcoh}} (\D ^\dag _{\PP,\Q})$,
la sous-catégorie pleine de
($F$-)$D ^\mathrm{b} _{\mathrm{coh}} (\D ^\dag _{\PP, \Q})$ des complexes $\D ^\dag _{\PP,\Q}$-surcohérents.
\end{defi}

\begin{vide}
  [Stabilité de la $\D ^\dag _{\PP,\Q}$-surcohérence et $\D$-modules arithmétiques surcohérents]
Avec les notations de \ref{defi-311surcoh},
comme $f$ est lisse et $T'$ est un diviseur de $P'$, les foncteurs $f ^*$ et $(\hdag T')$ sont exacts.
Il en résulte qu'un complexe est $\D ^\dag _{\PP,\Q}$-surcohérent si et seulement si ses espaces de cohomologie le sont.

De plus, d'après respectivement \cite[3.1.5]{caro_surcoherent}, \cite[3.1.7]{caro_surcoherent}, \cite[3.1.9]{caro_surcoherent},
la $\D ^\dag _{\PP,\Q}$-surcohérence est stable par foncteur cohomologique local, image directe par un morphisme propre,
image inverse extraordinaire.

Donnons ici une application de cette stabilité de la surcohérence (voir \cite[3.2]{caro_surcoherent}).
Soit $Y$ une $k$-variété telle qu'il existe un plongement de $Y$ dans un $\V$-schéma formel propre et lisse  $\PP$,
un diviseur $T$ de $P$ tel que $Y = \overline{Y} \setminus T$, où $\overline{Y}$ est l'adhérence schématique
de $Y$ dans $P$. La stabilité de la surcohérence permet d'établir que la catégorie des
$\D ^\dag _{\PP,\Q}$-modules surcohérents tels que $\E \riso \E (\hdag T)$ et
$\R \underline{\Gamma} ^\dag _{\overline{Y}} (\E) \riso \E$ ne dépend que de $Y$ (à isomorphisme canonique près), i.e., ne dépende pas du choix
de tels plongements de $Y$ dans un $\V$-schéma formel propre et lisse.
Ses objets sont alors naturellement appelés {\og $\D _Y$-modules arithmétiques surcohérents\fg}.
Plus généralement, lorsque $Y$ est une $k$-variété quelconque, on construit par recollement
cette catégorie des $\D _Y$-modules arithmétiques surcohérents.

\end{vide}

Il est possible de définir une notion de la surcohérence en évitant de recourir au foncteur cohomologique local de la façon suivante :
\begin{defi}
\label{defi-SC}
    Soient
    $X$ un sous-schéma fermé de $P$,
$\E \in (F\text{-})D ( \D ^\dag _{\PP,\Q})$.
On dit que $\E$ {\og est surcohérent sur $X$\fg} si
$\E \in D ^\mathrm{b} _\mathrm{coh} ( \D ^\dag _{\PP,\Q})$ et
pour tout morphisme $f$ : $\PP' \rightarrow \PP$ de $\V$-schémas formels lisses tel que
$f (P') \subset X$,
    $f ^! (\E) \in D ^\mathrm{b} _\mathrm{coh} ( \D ^\dag _{\PP',\Q})$.
\end{defi}

\begin{exem}
\label{exem}
\begin{enumerate}
  \item \label{exem-i}
  Soit $\E $ est un isocristal convergent sur $\PP$, i.e., un $\D ^\dag _{\PP,\Q}$-module cohérent, $\O _{\PP,\Q}$-cohérent.
  Il est facile de voir que $\E$ est surcohérent sur $P$ (a priori, ne pas confondre avec {\og surcohérent sur $\PP$\fg}).
  En effet, pour tout morphisme $f$ : $\PP' \rightarrow \PP$ de $\V$-schémas formels lisses,
  $f ^! (\E) \riso f ^* (\E) [d_{P'/P}]$, où $f ^* (\E) $ est l'isocristal convergent sur $\PP'$ déduit de $\E$
  par image inverse par $f$.

  De même, si $\E \in D ^\mathrm{b} _\mathrm{coh} ( \D ^\dag _{\PP,\Q}) \cap D ^\mathrm{b} _\mathrm{coh} ( \O _{\PP,\Q})$, alors $\E$ est surcohérent sur
  $P$.

  \item \label{exem-ii} Soit $\E\in D ^\mathrm{b} _\mathrm{surcoh} ( \D ^\dag _{\PP,\Q})$. Alors, pour tout sous-schéma fermé $X$ de $P$,
  $\E$ est surcohérent sur $X$. En effet, cela résulte de la stabilité de la surcohérence
par image inverse extraordinaire (\cite[3.1.7]{caro_surcoherent}).
Lorsque $\PP$ est de surcroît propre et $\E$ est muni d'une structure de Frobenius,
on vérifiera via \ref{carac-surcoh} que la réciproque est exacte.

\end{enumerate}

\end{exem}

\begin{rema}
\begin{enumerate}

\item  Avec les notations \ref{defi-SC},
  il n'est pas immédiat qu'un complexe $\E \in F\text{-}D ^\mathrm{b} _\mathrm{coh} ( \D ^\dag _{\PP,\Q})$
  soit surcohérent sur $X$ si et seulement si, pour tout entier $j$,
  $\H ^j (\E)$ soit surcohérent sur $X$.
  Le problème vient du fait que, comme le morphisme $f$ n'est pas forcément lisse,
  le foncteur $\H ^j$ ne commute pas en général à $f ^!$.
  Pour définir la {\og $\D ^\dag _{\PP,\Q}$-surcohérence\fg},
  nous avions pris soin de construire une notion stable par les foncteurs $\H ^j$.
  Par contre, d'après \ref{carac-surcoh} c'est bien le cas pour la surcohérence sur $P$ lorsque
  $\PP$ est propre.

\item Si $\E$ est surcohérent sur $X$, il n'est pas immédiat que
$\R \underline{\Gamma} _{X} ^\dag (  \E)$ soit surcohérent sur $X$.
D'après le lemme \ref{carac-surcoh-prelem1},
tout le problème est de vérifier
que $\R \underline{\Gamma} _{X} ^\dag (  \E)\in  D ^\mathrm{b} _\mathrm{coh} ( \D ^\dag _{\PP,\Q})$.

\end{enumerate}
\end{rema}
\medskip

Le lemme qui suit est immédiat.

\begin{lemm}
\label{lemm0}
\begin{enumerate}
  \item \label{lemm0i}
Soient $X$ un sous-schéma fermé de $P$, $\E \in D ( \D ^\dag _{\PP,\Q})$.
Pour tout recouvrement ouvert $(\PP _\alpha ) _{\alpha \in \Lambda}$ de $\PP$,
$\E$ est surcohérent sur $X$
si et seulement si, pour tout $\alpha \in \Lambda$,  $\E | \PP _\alpha$ est surcohérent sur $X \cap P _\alpha$.

\item \label{lemm0ii}
Soient $X$ un sous-schéma fermé de $P$ et
$\E' \rightarrow \E \rightarrow \E'' \rightarrow \E'[1]$ un triangle distingué de
$D ( \D ^\dag _{\PP,\Q})$.
Si deux des complexes sont surcohérents sur $X$, alors il en
est de même du troisième.

\end{enumerate}

\end{lemm}

\begin{lemm}
\label{carac-surcoh-lemele}
Soient $f$ : $\PP ' \rightarrow \PP$ un morphisme lisse de $\V$-schémas formels lisses,
$X$, $X'$ deux sous-schémas fermés de respectivement $P$ et $P'$ tels que $f$ induise
un morphisme $X' \rightarrow X$.

Si $\E $ est un complexe de $D ^\mathrm{b} _\mathrm{coh} ( \D ^\dag _{\PP,\Q})$
surcohérent sur $X$ alors $f ^! (\E) $ est un complexe
surcohérent sur $X'$.
  \end{lemm}
\begin{proof}
  Comme $f$ est lisse, $f ^! (\E) \in D ^\mathrm{b} _\mathrm{coh} ( \D ^\dag _{\PP',\Q})$.
  Le fait que $f ^! (\E) $ soit surcohérent sur $X'$ résulte alors de la transitivité pour la composition de l'image inverse extraordinaire.
\end{proof}

\begin{lemm}
  \label{carac-surcoh-prelem1}
      Soient
      $X$ un sous-schéma fermé de $P$,
$\E \in D ( \D ^\dag _{\PP,\Q})$.

Alors, $\R \underline{\Gamma} _{X} ^\dag (  \E)$ est surcohérent sur $X$ si et seulement si
$\R \underline{\Gamma} _{X} ^\dag (  \E)\in D ^\mathrm{b} _\mathrm{coh} ( \D ^\dag _{\PP,\Q})$
et $\E $ est surcohérent sur $X$.
\end{lemm}

\begin{proof}
  Soit $f$ : $\PP' \rightarrow \PP$ un morphisme de $\V$-schémas formels lisses tel que
$f (P') \subset X$. Le lemme résulte de l'isomorphisme
$f ^! \circ \R \underline{\Gamma} _{X} ^\dag (  \E)
\riso
\R \underline{\Gamma} _{f^{-1}(X)} ^\dag \circ f ^! (  \E)
=\R \underline{\Gamma} _{P'} ^\dag \circ f ^! (  \E)
=f ^! (  \E)$.
\end{proof}

Dans le cas où le sous-schéma fermé est lisse, le lemme \ref{carac-surcoh-prelem1} précédent s'affine de la manière suivante.
\begin{lemm}
\label{carac-surcoh-lem1}
  Soient
  $X$ un sous-schéma fermé lisse de $P$,
$\E \in D ^\mathrm{b} _\mathrm{coh} ( \D ^\dag _{\PP,\Q})$.
\begin{enumerate}
\item  \label{carac-surcoh-lem1-1}   Si $X _1,\dots, X _n$ sont les composantes irréductibles de $X$,
$\E $ est surcohérent sur $X$ si et seulement si, pour tout $i =1,\dots , n$,
    $\E$ est surcohérent sur $X _i$.

\item \label{carac-surcoh-lem1-2} On suppose que l'immersion fermée $X \hookrightarrow P$ se relève en
un morphisme $u\,:\,\X \hookrightarrow \PP$ de $\V$-schémas formels lisses.
Alors,
$\E$ est surcohérent sur $X$ si et seulement si $u ^! (\E)$ est surcohérent sur $X$.

  \item \label{carac-surcoh-lem1-3} Le complexe $\E$ est surcohérent sur $X$ si et seulement si
  $\R \underline{\Gamma} _{X} ^\dag (  \E)$
  est surcohérent sur $X$.

  \item \label{carac-surcoh-lem1-4} Si $\E$ est surcohérent sur $X$ alors,
  pour tout diviseur à croisements normaux strict $X'$ de $X$,
$  \R \underline{\Gamma} _{X'} ^\dag (  \E) $
est surcohérent sur $X'$.
\end{enumerate}
\end{lemm}
\begin{proof}
La première assertion provient du fait que la surcohérence sur $X$ est une notion locale (voir \ref{lemm0}.\ref{lemm0i}). La seconde est immédiate.
\'Etablissons à présent la troisième assertion.
Par \ref{carac-surcoh-prelem1}, si $\R \underline{\Gamma} _{X} ^\dag (  \E)$ est surcohérent sur $X$ alors
$\E $ est surcohérent sur $X$. Réciproquement, supposons $\E$ surcohérent sur $X$.
Comme le fait que $\R \underline{\Gamma} _{X} ^\dag (  \E)$ soit surcohérent sur $X$ est local en $\PP$ (voir \ref{lemm0}.\ref{lemm0i}),
on peut supposer que l'immersion fermée $X \hookrightarrow P$ se relève
en un morphisme $u$ : $\X \hookrightarrow \PP$ de $\V$-schémas formels lisses.
Or, d'après \ref{carac-surcoh-lem1}.\ref{carac-surcoh-lem1-2},
$u^! (\E)$ est cohérent. Il en résulte que
$\R \underline{\Gamma} _{X} ^\dag (  \E) \riso u _+ u^! (\E)$ est aussi cohérent.
Par \ref{carac-surcoh-prelem1}, cela implique que
$\R \underline{\Gamma} _{X} ^\dag (  \E)$ est surcohérent sur $X$.

Traitons à présent la dernière assertion.
On suppose donc que $\E$ est surcohérent sur $X$.
Comme $\E$ est a fortiori surcohérent sur $X'$ (voir \ref{carac-surcoh-lemele}),
d'après \ref{carac-surcoh-prelem1}, il s'agit d'établir que
$  \R \underline{\Gamma} _{X'} ^\dag (  \E) \in   D ^\mathrm{b} _\mathrm{coh} ( \D ^\dag _{\PP,\Q})$.
On procède pour cela par récurrence sur le nombre de composantes
irréductibles $X'_1,\dots ,X ' _r$ de $X'$.
Supposons $r=1$. Comme $\E$ est surcohérent sur $X'$, d'après
\ref{carac-surcoh-lem1}.\ref{carac-surcoh-lem1-3},
$  \R \underline{\Gamma} _{X'} ^\dag (  \E) \in   D ^\mathrm{b} _\mathrm{coh} ( \D ^\dag _{\PP,\Q})$.
Supposons donc $r \geq 2$. Posons $X '' _1= \cup _{i=2,\dots, r} X' _i$.
On dispose alors du triangle distingué de Mayer-Vietoris (voir \cite[2.2.16.1]{caro_surcoherent}) :
\begin{equation}\label{MVlemmiv}
  \R \underline{\Gamma} ^\dag _{X' _1\cap X'' _1}(\E ) \rightarrow
  \R \underline{\Gamma} ^\dag _{X' _1}(\E ) \oplus
\R \underline{\Gamma} ^\dag _{X'' _1 }(\E )  \rightarrow
\R \underline{\Gamma} ^\dag _{X' _1\cup X'' _1}(\E ) \rightarrow
\R \underline{\Gamma} ^\dag _{X' _1 \cap X'' _1}(\E )[1].
\end{equation}
Par hypothèse de récurrence, $\R \underline{\Gamma} ^\dag _{X' _1}(\E ) \oplus
\R \underline{\Gamma} ^\dag _{X'' _1 }(\E ) $ est cohérent.
Comme
$X' _1\cap X'' _1$ est de plus un diviseur à croisements normaux de $X'_1$ à $r -1$ composantes irréductibles,
par hypothèse de récurrence,
$\R \underline{\Gamma} ^\dag _{X' _1\cap X'' _1}(\E ) $ est cohérent.
Il résulte du triangle distingué \ref{MVlemmiv} que
$\R \underline{\Gamma} ^\dag _{X' _1\cup X'' _1}(\E )$ est cohérent.
D'où le résultat.

\end{proof}

\subsection{Cas des complexes surcohérents sur les points fermés}

Soient $\PP$ un $\V$-schéma formel lisse, $X$ un sous-schéma fermé de $P$.

\begin{defi}
  Un complexe $\E\in D ^\mathrm{b} _\mathrm{coh} ( \D ^\dag _{\PP,\Q})$
  est dit {\og à fibres extraordinaires finies sur $X$\fg}
  si, pour tout point fermé $x$ de $X$, pour tout relèvement $i _x$ de l'immersion fermée induite par $x$ de la forme
  $i_x\,:\,\Spf \V (x) \hookrightarrow \PP$,
  les espaces de cohomologie de $i ^! _x (\E )$ sont des $K$-espaces vectoriels de dimension finie.

Lorsque $\E$ est un complexe à fibres extraordinaires finies sur $P$,
on dira simplement que $\E$ est {\og à fibres extraordinaires finies\fg}.

On remarque que si $\E$ est à fibres extraordinaires finies sur $X$, il n'est pas évident qu'il en est de même pour ses espaces de cohomologie.
\end{defi}

\begin{prop}
\label{eq-fextfini-surcohpt}
  Soit $\E\in D ^\mathrm{b} _\mathrm{coh} ( \D ^\dag _{\PP,\Q})$. Alors $\E$ est à fibres extraordinaires finies sur $X$
  si et seulement si $\E$ est surcohérent sur chacun des points fermés de $X$.
\end{prop}
\begin{proof}
  Soit $x$ un point fermé de $X$. Alors, d'après \ref{carac-surcoh-lem1}.\ref{carac-surcoh-lem1-2},
  $\E$ est surcohérent sur $\{x\}$ si et seulement si
  $i _x ^! (\E)$ est surcohérent sur $\{x\}$, ce qui équivaut à dire que les espaces de cohomologie de $i _x ^! (\E)$
  soient des $K$-espaces vectoriels de dimension finie (e.g., voir \ref{exem}.\ref{exem-i}).
\end{proof}
\medskip

Le théorème qui suit est une version légèrement plus fine que \cite[2.2.17]{caro_courbe-nouveau}.
\begin{theo}
  \label{2217bis}
Soit $\E$ un $F$-$\D ^{\dag } _{\PP , \Q }$-module cohérent tel que, pour tout point fermé $x$ de $P$,
le $K$-espace vectoriel $ i ^* _x ( \E ) $ est de dimension finie.
Il existe alors un diviseur $T$ de $P$ tel que $(\hdag T) (\E)$ soit un $F$-isocristal sur $P\setminus T$ surconvergent le long de $T$.
\end{theo}
\begin{proof}
Il suffit de reprendre la preuve du théorème \cite[2.2.17]{caro_courbe-nouveau}
en remplaçant {\og $i ^! _x$\fg} par {\og $i ^* _x$\fg}. En effet, seule la finitude de $i ^* _x (\E)$ est utilisée.

\end{proof}

\begin{coro}
  \label{fibextfinidense}
Soit $\E$ un complexe de $F\text{-}D ^\mathrm{b} _\mathrm{coh} ( \D ^\dag _{\PP,\Q})$
à support dans $X$ et surcohérent sur ses points fermés.
Il existe alors un ouvert $\U$ de $\PP$ tel que $U \cap X$ soit lisse et dense dans $X$ et tel que
les espaces de cohomologie de $\E |\U$
soient dans l'image essentielle du foncteur $\sp _{Y \hookrightarrow \U,+}$ (voir \ref{Fisocsurcoh}).
\end{coro}

\begin{proof}
Comme il existe un ouvert affine $\U$ de $\PP$ tel que $U \cap X$ soit lisse et dense dans $X$,
quitte à remplacer $\PP$ par $\U$, on se ramène au cas où $X$ est lisse
et $X\subset P$ se relève en un morphisme $u\,:\,\X \hookrightarrow \PP$ de $\V$-schémas formels lisses.
Dans ce cas, $\sp _{X \hookrightarrow \PP,+} = u_+ \circ \sp _{X \hookrightarrow \X,+}$ (on rappelle aussi que $\sp _* =\sp _{X \hookrightarrow \X,+}$,
où $\sp$ désigne le morphisme de spécialisation de la fibre générique de Raynaud de $\X$ vers $\X$ : voir \ref{Fisocsurcoh}).
D'après l'analogue $p$-adique de Berthelot du théorème de Kashiwara (voir \cite[5.3.3]{Beintro2}),
$u ^! \circ u _+ \riso Id$.
On peut donc supposer $X$ égal à $P$. Il ne coûte rien non
plus de supposer $P$ intègre. \'Etablissons à présent une récurrence sur le cardinal de l'ensemble
$\{j, \H ^j (\E) \not = 0\}$. Soit $N$ l'élément maximal de cet ensemble. Si le support de $\H ^N (\E) $ n'est pas égal à $P$, il
est alors de dimension strictement inférieure à $P$. Dans ce cas, en enlevant le support de $\H ^N (\E) $ à $P$,
l'hypothèse de récurrence nous permet de conclure. Traitons à présent le cas où le support de $\H ^N (\E) $ est égal à $P$.
Comme pour tout point fermé $x$ de $X$, pour tout entier $i \not \in \{0,\dots, d _P\}$, pour tout entier $j $, $\H ^i i_x ^! (\H ^j (\E)) =0$,
la deuxième suite spectrale d'hypercohomologie donne  $\H ^{d_P +N} i_x ^! (\E) \riso i ^* _x (\H ^N (\E) )$.
Ainsi, $i ^* _x (\H ^N (\E) )$ est un $K$-espace vectoriel de dimension finie.
D'après \ref{2217bis}, on en déduit qu'il existe un ouvert dense $\U$ de $\PP$
tel que $\H ^N (\E) |\U$ soit $\O _{\U,\Q}$-cohérent. Quitte à remplacer $\PP$ par $\U$, on peut donc supposer
que $\H ^N (\E) $ est $\O _{\PP,\Q}$-cohérent. Dans ce cas, pour tout $j \not =0$, $\H ^j \L i ^* _x (\H ^N (\E) )=0$.
On termine alors par dévissage (via le triangle distingué de troncation) grâce à l'hypothèse de récurrence.

\end{proof}

Lorsque $X$ est intègre et lisse, le corollaire \ref{fibextfinidense} s'affine de la manière suivante :
\begin{coro}
  \label{fibextfinidense2}
On suppose $X$ intègre et lisse.
Soit $\E$ un complexe de $F\text{-}D ^\mathrm{b} _\mathrm{coh} ( \D ^\dag _{\PP,\Q})$
à support dans $X$ et surcohérent sur ses points fermés.
Il existe alors un diviseur $T$ de $P$ tel que $X \setminus T$ soit dense dans $X$ et tel que
les espaces de cohomologie de $\E (\hdag T) $
soient dans l'image essentielle du foncteur $\sp _{X \hookrightarrow \PP,T,+}$ (voir \ref{Fisocsurcoh}).
\end{coro}

\begin{proof}
D'après \ref{fibextfinidense}, il existe un ouvert
$\U$ de $\PP$ tel que $U \cap X$ soit lisse et dense dans $X$ et tel que
les espaces de cohomologie de $\E |\U$
soient dans l'image essentielle du foncteur $\sp _{Y \hookrightarrow \U,+}$.
Quitte à rétrécir $U$, on peut en outre supposer que $P \setminus U$ est le support d'un diviseur $T$ (on remarque alors que
$T \cap X$ est aussi un diviseur de $X$).
Si $X \subset P$ se relève en un morphisme $u\,:\, \X \hookrightarrow \PP$ de
$\V$-schémas formels lisses, en notant $\Y$ l'ouvert de $\X$ complémentaire de $T \cap X$,
cela signifie que
les espaces de cohomologie de $u^! (\E (\hdag T)) |\Y$ sont
$\O _{\Y,\Q}$-cohérents.
D'après la caractérisation de Berthelot des isocristaux surconvergents
de \cite{Be4} (plus précisément, voir \cite[2.2.12]{caro_courbe-nouveau}),
cela implique que
les espaces de cohomologie de $u^! (\E (\hdag T))$ sont
$\O _{\X} (\hdag T \cap X) _{\Q}$-cohérents.
On conclut alors via la description de
l'image essentielle du foncteur $\sp _{X \hookrightarrow \PP,T,+}$
de \cite[2.5.10]{caro-construction}.
\end{proof}

\section{Sur l'équivalence entre surcohérence et surcohérence sur un fermé}

\subsection{Cas des isocristaux surconvergents}

\begin{vide}
  [Rappels sur les $\D$-modules arithmétiques associés aux isocristaux surconvergents et notations]
\label{Fisocsurcoh}
Soient $\PP$ un $\V$-schéma formel séparé et lisse, $T$ un diviseur de $P$, $\U := \PP \setminus T$,
$X$ un sous-schéma fermé. On suppose que $Y:=X \setminus T$ est lisse.

$\bullet$ On désigne par $(F\text{-})\mathrm{Isoc} ^{\dag}(Y,X/K )$, la catégorie
des ($F$-)isocristaux sur $Y$ surconvergents le long de $T \cap X$.
Lorsque $X$ est propre, la catégorie $(F\text{-})\mathrm{Isoc} ^{\dag}(Y,X/K )$
est indépendante du choix de la compactification $X$ de $Y$ et se note
$(F\text{-})\mathrm{Isoc} ^{\dag}(Y/K )$. Ses objets sont appelés
($F$-)isocristaux sur $Y$.

$\bullet$ Lorsque $X$ est lisse et $T \cap X$ est un diviseur de $X$, on dispose d'un foncteur
pleinement fidèle noté $\sp _{X \hookrightarrow \PP,T,+}$ de la catégorie
$(F\text{-})\mathrm{Isoc} ^{\dag}(Y,X/K )$
dans celle des $(F\text{-})\D ^\dag _{\PP} (\hdag T) _{\Q}$-modules cohérents à support dans $X$.
En gros, le foncteur $\sp _{X \hookrightarrow \PP,T,+}$ se construit de la manière suivante :
lorsque $X \hookrightarrow P$ se relève en un morphisme de $\V$-schémas formels lisses de la forme
$u\,:\, \X \hookrightarrow \PP$ (e.g., si $\PP$ est affine),
alors $\sp _{X \hookrightarrow \PP,T,+}= u_+ \circ \sp _*$, où $\sp$ désigne le morphisme de spécialisation
de la fibre générique de Raynaud de $\PP$ vers $\PP$.
En vérifiant (grâce aux propriétés d'indépendance des isocristaux surconvergents et
des $\D$-modules arithmétiques) que ce foncteur ne dépend canoniquement pas du choix du relèvement $u$,
on obtient $\sp _{X \hookrightarrow \PP,T,+}$ par recollement.

$\bullet$
On note $\DD _T$ le foncteur dual $\smash{\D} ^\dag _{\PP} (\hdag T) _\Q$-linéaire.
La catégorie $F$-$\mathrm{Isoc} ^{\dag \dag}( \PP, T, X/K)$
est alors définie de la manière suivante
\cite[6.2.1]{caro_devissge_surcoh}
  \begin{itemize}
    \item les objets sont les
  $F$-$\smash{\D} ^\dag _{\PP} (\hdag T) _\Q$-modules surcohérents $\E$ à support dans $X$ tels que
  $\DD _T (\E)$ soit $\smash{\D} ^\dag _{\PP} (\hdag T) _\Q$-surcohérent et
  tels qu'il existe un isocristal convergent $G$ sur $Y$
vérifiant
    $\E |{\U} \riso \sp _{Y\hookrightarrow \U, +} (G)$ ;
   \item les flèches sont les morphismes $F$-$\smash{\D} ^\dag _{\PP} (\hdag T) _\Q$-linéaires.
  \end{itemize}

Lorsque $\PP$ est propre,
cette catégorie $F$-$\mathrm{Isoc} ^{\dag \dag}( \PP, T, X/K)$ ne dépend canoniquement pas du choix
de $(\PP, T,X)$ tels que $X \setminus T=Y$. Elle sera alors notée
$F$-$\mathrm{Isoc} ^{\dag \dag}( Y/K)$. Ses objets sont les {\og $F$-isocristaux surcohérents sur $Y$\fg}.
\medskip

$\bullet$ Lorsque $Y$ est un $k$-schéma séparé lisse quelconque, on construit plus généralement la catégorie
$F$-$\mathrm{Isoc} ^{\dag \dag}( Y/K)$ par recollement (voir \cite[2.2.4]{caro-2006-surcoh-surcv}).
De plus, la catégorie ($F$-)$\mathrm{Isoc} ^{\dag}( Y/K)$ des ($F$-)isocristaux surconvergents sur $Y$
de Berthelot est encore définie.

\end{vide}

$\bullet$ Soit $Y $ un $k$-schéma séparé et lisse.
D'après \cite[2.3.1]{caro-2006-surcoh-surcv},
on dispose de l'équivalence canonique de catégories
\begin{equation}
  \label{eqcat-gen}
\sp _{Y,+}\ :\ F\text{-}\mathrm{Isoc} ^{\dag}(Y/K )\cong F\text{-}\mathrm{Isoc} ^{\dag \dag}(Y/K ).
\end{equation}

La proposition suivante
nous sera utile afin d'établir la première étape (i.e., $(Q _{n-1}) \Rightarrow (P _n)$) de la preuve de \ref{carac-surcoh}.
Remarquons d'ailleurs qu'ici la structure de Frobenius n'est pas indispensable.
\begin{prop}
\label{carac-surcoh-lem2}
Soient $\PP$ un $\V$-schéma formel séparé et lisse, $X$ un sous-schéma fermé lisse de $P$,
  $T$ un diviseur de $P$ tel que $T \cap X$ soit un diviseur de $X$.
Soit $\E $ un complexe de $D ^\mathrm{b} _\mathrm{coh} ( \D ^\dag _{\PP} (\hdag T) _{\Q})$
tel que,
pour tout entier $j$, $\H ^j (\E)$ soit dans l'image essentielle de $\sp _{X \hookrightarrow \PP, T,+}$.

Alors, le complexe $\E$ est surcohérent sur $X$ si et seulement si $\E \in D ^\mathrm{b} _\mathrm{surcoh} ( \D ^\dag _{\PP,\Q})$.
\end{prop}

\begin{proof}
Il résulte de \ref{exem}.\ref{exem-ii} (i.e., \cite[3.1.7]{caro_surcoherent})
 que si $\E \in D ^\mathrm{b} _\mathrm{surcoh} ( \D ^\dag _{\PP,\Q})$ alors
$\E$ est surcohérent sur $X$. Réciproquement, supposons $\E$ surcohérent sur $X$.
Comme l'assertion est locale en $\PP$, on se ramène au cas où $\PP$ est affine.
L'immersion fermée $X \hookrightarrow P$ se relève donc en
un morphisme $u\,:\,\X \hookrightarrow \PP$ de $\V$-schémas formels lisses.
D'après \ref{carac-surcoh-lem1}.\ref{carac-surcoh-lem1-2},
$ u^! (\E)$ est surcohérent sur $X$.
Comme $\E$ est à support dans $X$, d'après le théorème de Kashiwara
$\E \riso u_+ u^! (\E)$. Comme la surcohérence est stable par le foncteur $u_+$,
on se ramène ainsi au cas où $X =P$.
Par \ref{carac-surcoh-lem1}.\ref{carac-surcoh-lem1-1}, on peut en outre supposer $X$ intègre.
On note dans ce cas $\X$ pour $\PP$.

Soit $g\,\: \, \widetilde{\X} \rightarrow \X$ un morphisme lisse de $\V$-schémas formels lisses.
D'après \ref{carac-surcoh-lemele},
$g ^! (\E)$ est surcohérent sur $\widetilde{X}$.
De plus, en notant $\widetilde{T} :=g ^{-1} (T)$,
alors $g ^! (\E)\riso (\hdag \widetilde{T})(g ^! (\E))$ et ses espaces de cohomologie
sont dans l'image essentielle du foncteur $\sp _* = \sp _{\widetilde{X} \hookrightarrow \widetilde{\X}, \widetilde{T},+}$
(par commutation des foncteurs de la forme $\sp _{X \hookrightarrow \PP, T,+}$ aux images inverses : voir \cite[4.1.8]{caro-construction}).
Pour prouver la surcohérence de $\E$,
on se ramène alors au cas où $g =Id$, i.e.,
il suffit de prouver que pour tout diviseur $Z$ de $X$,
$\E (\hdag Z) \in D ^\mathrm{b} _\mathrm{coh} ( \D ^\dag _{\X,\Q})$.
Comme $\E (\hdag Z)  \riso \E (\hdag Z \cup T) $, on peut en outre supposer $T \subset Z$.

  Grâce au théorème de désingularisation de de Jong (voir \cite{dejong}),
  il existe un diagramme commutatif de la forme
\begin{equation}
    \label{diag631dejongpre}
  \xymatrix @R=0,3cm {
  { Y '} \ar[r] \ar[d] ^b & {X'} \ar[r] ^{u'} \ar[d] ^a & {\PP'} \ar[d] ^f \\
  {Y} \ar[r] & {X } \ar[r] ^{u} & {\X,}}
\end{equation}
  où $f$ est un morphisme propre et lisse de $\V$-schémas formels lisses, $Y := X \setminus Z$,
  le carré de gauche est cartésien, $X'$ est un $k$-schéma lisse, $u'$ est une immersion fermée,
  $a$ est un morphisme
  projectif, surjectif, génériquement fini et étale tel que
  $a ^{-1} (Z)$ soit un diviseur à croisement normaux strict de $X'$.
De plus, comme pour tout entier $j$, $\H ^j \E (\hdag Z)$ est dans l'image essentielle de
$\sp _* = \sp _{X\hookrightarrow \X, Z,+}$
(i.e., est associé à un isocristal sur $Y$ surconvergent le long de $Z$),
d'après la preuve de \cite[6.1.4]{caro_devissge_surcoh},
le faisceau $\H ^j \E (\hdag Z)$ est alors un facteur direct de
$f _+  \R \underline{\Gamma} _{X'} ^\dag  f ^! (\H ^j \E (\hdag Z) )$.
Il suffit donc d'établir que $f _+  \R \underline{\Gamma} _{X'} ^\dag  f ^! \H ^j (\E (\hdag Z) )
\in
D ^\mathrm{b} _\mathrm{coh} ( \D ^\dag _{\X,\Q})$.

Le triangle distingué de localisation de $\R \underline{\Gamma} _{X'} ^\dag ( f ^! \E)$
en $Z':=f ^{-1}(Z)$ s'écrit (on utilise aussi \cite[2.2.18]{caro_surcoherent}) :
\begin{equation}
\label{tri-locZ'}
  \R \underline{\Gamma} _{Z '\cap X'} ^\dag ( f ^! \E)
  \rightarrow
  \R \underline{\Gamma} _{X'} ^\dag ( f ^! \E)
  \rightarrow
\R \underline{\Gamma} _{X'} ^\dag  f ^! (\E (\hdag Z) )
  \rightarrow
  \R \underline{\Gamma} _{Z '\cap X'} ^\dag ( f ^! \E) [1].
\end{equation}
Par \ref{carac-surcoh-lemele}, $f ^! \E$ est surcohérent sur $X'$.
Il résulte alors respectivement de \ref{carac-surcoh-lem1}.\ref{carac-surcoh-lem1-3}
et \ref{carac-surcoh-lem1}.\ref{carac-surcoh-lem1-4} que
$\R \underline{\Gamma} _{X'} ^\dag ( f ^! \E), \R \underline{\Gamma} _{Z '\cap X'} ^\dag ( f ^! \E)\in
D ^\mathrm{b} _\mathrm{coh} ( \D ^\dag _{\PP',\Q})$.
Via \ref{tri-locZ'}, il en découle que
$\R \underline{\Gamma} _{X'} ^\dag  f ^! (\E (\hdag Z) )\in
D ^\mathrm{b} _\mathrm{coh} ( \D ^\dag _{\PP',\Q})$.
Puisque, pour tout entier $j$, $\H ^j  \E (\hdag Z)$ est un isocristal sur $Y$ surconvergent le long de $Z$,
par \cite[4.1.8]{caro-construction}, pour tout entier $i \not =0$,
$\H ^i \R \underline{\Gamma} _{X'} ^\dag  f ^! (\H ^j  \E (\hdag Z) )=0$.
Via la deuxième suite spectrale d'hypercohomologie de
$\R \underline{\Gamma} _{X'} ^\dag  f ^!$ appliquée à $\E (\hdag Z)$, il en résulte
$\H ^j (\R \underline{\Gamma} _{X'} ^\dag  f ^! (\E (\hdag Z) ))
\riso
\R \underline{\Gamma} _{X'} ^\dag  f ^! (\H ^j  \E (\hdag Z) )$.
Ainsi, $\R \underline{\Gamma} _{X'} ^\dag  f ^! (\H ^j  \E (\hdag Z) ) $ est
un $\D ^\dag _{\PP',\Q}$-module cohérent.
Par stabilité de la cohérence par image directe par un morphisme propre, on en déduit que
$f _+  \R \underline{\Gamma} _{X'} ^\dag  f ^! \H ^j (\E (\hdag Z) )
\in
D ^\mathrm{b} _\mathrm{coh} ( \D ^\dag _{\X,\Q})$.
D'où le résultat.

\end{proof}

\subsection{Cas général}
\begin{theo}
\label{carac-surcoh}
 Soient $\PP$ un $\V$-schéma formel propre et lisse, $X$ un sous-schéma fermé de $P$,
$\E \in F\text{-}D ^\mathrm{b} _\mathrm{coh} ( \D ^\dag _{\PP,\Q})$ à support dans $X$.
  Les deux assertions suivantes sont alors équivalentes :
  \begin{enumerate}
    \item \label{carac-surcoh1} Le $F$-complexe $\E$ est surcohérent sur $X$.
    \item \label{carac-surcoh2} Le $F$-complexe $\E$ appartient à $ F\text{-}D ^\mathrm{b} _\mathrm{surcoh} ( \D ^\dag _{\PP,\Q})$.
  \end{enumerate}
\end{theo}

\begin{proof}
Il résulte de \ref{exem}.\ref{exem-ii} (i.e., \cite[3.1.7]{caro_surcoherent})
 l'implication $(\ref{carac-surcoh2})\Rightarrow (\ref{carac-surcoh1})$.
Traitons à présent la réciproque. Considérons pour tout entier
$n$ les deux propriétés suivantes ($\E$, $\PP$ et $X$ pouvant naturellement variés mais vérifiant les hypothèses du théorème) :
\begin{enumerate}
\item [$(P _n)$] Si $X$ est lisse, $\dim X \leq n$ et $\E$ est surcohérent sur $X$
  alors $\E \in F\text{-}D ^\mathrm{b} _\mathrm{surcoh} ( \D ^\dag _{\PP,\Q})$.

  \item [$(Q _n)$]
  Si $\dim X \leq n$ et $\E$ est surcohérent sur $X$ alors
  $\E \in F\text{-}D ^\mathrm{b} _\mathrm{surcoh} ( \D ^\dag _{\PP,\Q})$.
\end{enumerate}

En gros, nous allons prouver l'implication $(Q _{n-1}) \Rightarrow (P _n)$
en nous ramenant par dévissage au cas des $F$-isocristaux surconvergents déjà traité dans \ref{carac-surcoh-lem2}.
Puis, on établit $\{(Q _{n-1})\text{ et }(P _n) \}\Rightarrow (Q _n)$ en désingularisant $X$ à la de Jong et
en faisant de la descente de $F$-isocristaux surconvergents.

$1 ^\circ)$ \'Etablissons maintenant l'implication $(Q _{n-1}) \Rightarrow (P _n)$.
Supposons ainsi que $X$ est lisse, $\dim X \leq n$ et $\E$ est surcohérent sur $X$.
Par \ref{carac-surcoh-lem1}.\ref{carac-surcoh-lem1-1}, il ne coûte rien de supposer $X$ intègre.
D'après \ref{fibextfinidense2},
comme $\E$ est à support dans $X$ et est à fibres extraordinaires finies, il existe
un diviseur $T $ de $P$ tel que $T \cap X$ soit un diviseur de $X$ et, pour tout
  entier $j$, $\H ^j (\hdag T) (\E)$ soit associé à un $F$-isocristal surconvergent sur $Y:=X \setminus T$,
  i.e., soit dans l'image essentielle du foncteur $\sp _{X \hookrightarrow \PP,T,+}$.
  De plus, il résulte de \cite[6.1.4]{caro_devissge_surcoh} que
  $\H ^j (\hdag T) (\E)\in F\text{-}\mathrm{Isoc} ^{\dag \dag}( \PP, T, X/K)$.
  Par contre, on ne peut directement utiliser \ref{carac-surcoh-lem2}, car
  il n'est pas immédiat que $(\hdag T) (\E)$ soit surcohérent sur $X$
  (plus précisément que $(\hdag T) (\E) \in F\text{-}D ^\mathrm{b} _\mathrm{coh} ( \D ^\dag _{\PP,\Q})$).
Heureusement, lorsque $T \cap X$ est un diviseur à croisements normaux strict de $X$,
les choses se simplifient grâce au lemme suivant :
\begin{lemm*}
\label{carac-surcoh-lem12}
Avec les notations de \ref{carac-surcoh}, supposons $X$ lisse et soit
  $X'$ un diviseur à croisements normaux strict de $X$.
  Si $\dim X \leq n$, $\E$ est surcohérent sur $X$
et $(Q _{n-1})$ est vraie, alors
  $\R \underline{\Gamma} _{ X'} ^\dag (  \E) \in
 F\text{-} D ^\mathrm{b} _\mathrm{surcoh} ( \D ^\dag _{\PP,\Q})$.
\end{lemm*}

\begin{proof}
Comme $X'$ est un diviseur à croisement normaux strict de $X$,
d'après \ref{carac-surcoh-lem1}.\ref{carac-surcoh-lem1-4},
$\R \underline{\Gamma} _{ X'} ^\dag (  \E) $ est surcohérent sur $X'$.
Comme $\dim X' \leq n-1$, l'hypothèse que $(Q _{n-1})$ soit vraie nous permet de conclure.

\end{proof}

L'idée est maintenant de se ramener au cas où $T \cap X$ est un diviseur à croisements normaux strict de $X$
(cas qui permet d'utiliser le lemme \ref{carac-surcoh-lem12}) de la façon suivante :
d'après le théorème de désingularisation de de Jong, il existe un diagramme commutatif de la forme
\begin{equation}
\notag
    \label{diag631dejong}
  \xymatrix @R=0,3cm {
  {X'} \ar[r] ^{u'} \ar[d] ^a & {\PP'} \ar[d] ^f \\
   {X } \ar[r] ^{u} & {\PP,}}
\end{equation}
  où $\PP'$ est un espace projectif formel de base $\PP$, $f$ est la projection canonique,
$X'$ est un $k$-schéma lisse, $u'$ est une immersion fermée,
  $a$ est un morphisme
  projectif, surjectif, génériquement fini et étale tel que
  $a ^{-1} (T\cap X )$ soit un diviseur à croisement normaux strict de $X'$.
  Notons $T' = f ^{-1} (T)$.
D'après \cite[6.3.1]{caro_devissge_surcoh},
comme $\H ^j \E (\hdag T) \in F\text{-}\mathrm{Isoc} ^{\dag \dag}( \PP, T, X/K)$,
$\H ^j \E (\hdag T)$ est un facteur direct de
$f _+  \R \underline{\Gamma} _{X'} ^\dag  f ^! (\H ^j \E (\hdag T) )$.

Par \ref{carac-surcoh-lemele}, $f ^! \E$ est surcohérent sur $X'$.
D'après \ref{carac-surcoh-lem1}.\ref{carac-surcoh-lem1-3},
$\R \underline{\Gamma} _{X'} ^\dag ( f ^! \E)$
est alors surcohérent sur $X'$.
De plus, comme $(Q _{n-1}) $ est vraie,
  il résulte du lemme \ref{carac-surcoh-lem12}
  que $\R \underline{\Gamma} _{T '\cap X'} ^\dag ( f ^! \E) \in
  F\text{-}D ^\mathrm{b} _\mathrm{surcoh} ( \D ^\dag _{\PP',\Q})$.
Or,
le triangle distingué de localisation (voir \ref{tri-loc-gen}) de $\R \underline{\Gamma} _{X'} ^\dag ( f ^! \E)$
en $T'$ s'écrit (on utilise aussi \cite[2.2.18]{caro_surcoherent}) :
\begin{equation}
\label{tri-locT'}
  \R \underline{\Gamma} _{T '\cap X'} ^\dag ( f ^! \E)
  \rightarrow   \R \underline{\Gamma} _{X'} ^\dag ( f ^! \E)
  \rightarrow
\R \underline{\Gamma} _{X'} ^\dag  f ^! (\E (\hdag T) )
  \rightarrow
  \R \underline{\Gamma} _{T '\cap X'} ^\dag ( f ^! \E) [1].
\end{equation}
On en déduit que
$\R \underline{\Gamma} _{X'} ^\dag  f ^! (\E (\hdag T) )$
est surcohérent sur $X'$ (voir \ref{lemm0}.\ref{lemm0ii}).

Comme $\H ^j (\E (\hdag T))$ est dans l'image essentielle de $\sp _{X \hookrightarrow \PP,T,+}$,
d'après \cite[4.1.8]{caro-construction},
$\R \underline{\Gamma} _{X'} ^\dag  f ^! (\H ^j  \E (\hdag T) )$ est
dans l'image essentielle du foncteur $\sp _{X '\hookrightarrow \PP',T',+}$.
En particulier,
pour tout entier $i \not =0$,
$\H ^i \R \underline{\Gamma} _{X'} ^\dag  f ^! (\H ^j  \E (\hdag T) )=0$.
Via la deuxième suite spectrale d'hypercohomologie de
$\R \underline{\Gamma} _{X'} ^\dag  f ^!$ appliquée à $\E (\hdag T)$, il en résulte
$\H ^j (\R \underline{\Gamma} _{X'} ^\dag  f ^! (\E (\hdag T) ))
\riso
\R \underline{\Gamma} _{X'} ^\dag  f ^! (\H ^j  \E (\hdag T) )$.
Ainsi,
$\R \underline{\Gamma} _{X'} ^\dag  f ^! (\E (\hdag T) )$
est un complexe surcohérent sur $X'$ dont les espaces de cohomologie
sont dans l'image essentielle de $\sp _{X '\hookrightarrow \PP',T',+}$.
D'après \ref{carac-surcoh-lem2},
il en résulte que
$\R \underline{\Gamma} _{X'} ^\dag  f ^! (\E (\hdag T) )
\in F\text{-}D ^\mathrm{b} _\mathrm{surcoh} ( \D ^\dag _{\PP,\Q})$.
On en déduit donc que
$\R \underline{\Gamma} _{X'} ^\dag  f ^! (\H ^j  \E (\hdag T) )$
est un
$F\text{-}\D ^\dag _{\PP',\Q}$-module surcohérent.
Comme la surcohérence est stable par image directe par un morphisme propre,
cela implique que
$f _+ \R \underline{\Gamma} _{X'} ^\dag  f ^! \H ^j (\E (\hdag T))\in
F\text{-}D ^\mathrm{b} _\mathrm{surcoh} ( \D ^\dag _{\PP,\Q})$.
Comme $\H ^j (\E (\hdag T))$ est un facteur direct de
$f _+ \R \underline{\Gamma} _{X'} ^\dag  f ^! \H ^j (\E (\hdag T))$,
$\H ^j (\E (\hdag T))$ est un $\D ^\dag _{\PP,\Q}$-module surcohérent.
Ainsi, $\E (\hdag T)\in
F\text{-}D ^\mathrm{b} _\mathrm{surcoh} ( \D ^\dag _{\PP,\Q})$.
Via le triangle de localisation de $\E$ en $T$ (voir \ref{tri-loc-gen}), il en dérive que
$\R \underline{\Gamma} _{T} ^\dag   (\E)$
est surcohérent sur $X$ (voir \ref{lemm0}.\ref{lemm0ii}).
Par \ref{carac-surcoh-lemele}, $\R \underline{\Gamma} _{T} ^\dag   (\E)$ est donc surcohérent sur $T \cap X$.
Comme en outre $\R \underline{\Gamma} _{T} ^\dag   (\E)$ est à support dans $T \cap X $,
comme $\dim T \cap X \leq n-1$, la validité de $(Q _{n-1})$ implique alors que
$\R \underline{\Gamma} _{T} ^\dag   (\E)\in
F\text{-}D ^\mathrm{b} _\mathrm{surcoh} ( \D ^\dag _{\PP,\Q})$.
On conclut alors via le triangle de localisation de $\E$ en $T$
que
$\E \in F\text{-}D ^\mathrm{b} _\mathrm{surcoh} ( \D ^\dag _{\PP,\Q})$.
\medskip

$2 ^\circ )$ Pour conclure, il suffit à présent de prouver l'implication $\{(Q _{n-1})\text{ et }(P _n) \}\Rightarrow (Q _n)$.
Supposons que $\dim X \leq n$ et $\E$ est surcohérent sur $X$.
Soit $Z $ une composante irréductible de $X$ de dimension égale à $\dim X$.
Il existe un diviseur $T$ de $P$ tel que $Z \setminus T=X \setminus T$
et $Y:=Z \setminus T$ soit affine, lisse et dense dans $Z$.
Par \ref{fibextfinidense}, comme $\E$ est surcohérent sur $Z$, quitte à augmenter $T$, on peut supposer
que, pour tout entier $j$, en notant $\U := \PP \setminus T$,
il existe un (unique à isomorphisme non canonique près)
$F$-isocristal convergent $G_j$ sur $Y$
tel que
$\H ^j \E |\U \riso \sp _{Y \hookrightarrow \U,+} (G_j)$.
On dispose grâce au théorème de désingularisation de de Jong du diagramme commutatif :
\begin{equation}
    \label{diag631dejong2}
  \xymatrix @R=0,3cm {
  { Y '} \ar[r] \ar[d] ^b & {Z'} \ar[r] ^{u'} \ar[d] ^a & {\PP'} \ar[d] ^f \\
  {Y} \ar[r] & {Z } \ar[r] ^{u} & {\PP,}}
\end{equation}
  où $\PP'$ est un espace projectif formel de base $\PP$, $f$ est la projection canonique,
  le carré de gauche est cartésien, $Z'$ est un $k$-schéma lisse, $u'$ est une immersion fermée,
  $a$ est un morphisme
  projectif, surjectif, génériquement fini et étale.
  Quitte à agrandir $T$, on peut en fait supposer que le morphisme $b$ est fini et étale.
Notons $T':=f ^{-1} (T)$, $\U' := \PP' \setminus T'$ et $g \,:\, \U' \rightarrow \U$ le morphisme induit
par $f$.

Comme $f$ est lisse et $\E(\hdag T)$ est $\D ^\dag _{\PP} (\hdag T) _{\Q}$-cohérent,
$f ^! (\E(\hdag T))$ est donc $\D ^\dag _{\PP'} (\hdag T') _{\Q}$-cohérent.
Comme $f$ est en outre propre, on obtient alors par adjonction (voir \ref{rapnot}) :
  $f _+ f ^! (\E(\hdag T)) \rightarrow \E(\hdag T)$. On en déduit :
  $f _+ \R \underline{\Gamma} ^\dag _{Z'} f ^! (\E (\hdag T)) \rightarrow \E (\hdag T)$.
En lui appliquant $\H ^j $, cela donne la flèche :
$\H ^j f _+  \R \underline{\Gamma} ^\dag _{Z'} f ^! (\E (\hdag T)) \rightarrow \H ^j  \E (\hdag T)$.
Nous allons maintenant prouvé que
$\H ^j f _+  \R \underline{\Gamma} ^\dag _{Z'} f ^! (\E (\hdag T))  \riso
f _+  (\H ^j  \R \underline{\Gamma} ^\dag _{Z'} f ^! (\E (\hdag T))) $.
Nous aurons pour cela besoin du lemme ci-après :
\begin{lemm*}
\label{carac-surcoh-lemm*2}
  Avec les notations de l'étape $2 ^\circ )$ de la preuve de \ref{carac-surcoh},
  pour tout entier $j$,
  $\H ^j \R \underline{\Gamma} ^\dag _{Z'} f ^! (\E (\hdag T))$
est un $\D ^\dag _{\PP',\Q}$-module surcohérent. De plus, il existe
 un $F$-isocristal $E'_j$ surconvergent sur $Y'$
tel que
$ \H ^j \R \underline{\Gamma} ^\dag _{Z'} f ^! (\E (\hdag T))\riso \sp _{Z'\hookrightarrow \PP', T',+} (E'_j )$
et
$\widehat{E'} _j \riso b ^* (G _j)$, où
$\widehat{E'} _j $ est le $F$-isocristal convergent sur $Y'$ induit par $E'_j$.
\end{lemm*}
\begin{proof}

  Comme $\E$ est surcohérent sur $Z$, par \ref{carac-surcoh-lemele} et \ref{carac-surcoh-lem1}.\ref{carac-surcoh-lem1-3},
$\R \underline{\Gamma} ^\dag _{Z'} f ^! (\E )$
est surcohérent sur $Z'$.
  L'hypothèse $(P _n)$ implique alors que $\R \underline{\Gamma} ^\dag _{Z'} f ^! (\E )
\in F\text{-}D ^\mathrm{b} _\mathrm{surcoh} ( \D ^\dag _{\PP',\Q})$.
Comme $\R \underline{\Gamma} ^\dag _{Z'} f ^! (\E (\hdag T))\riso
(\hdag T' ) \R \underline{\Gamma} ^\dag _{Z'} f ^! (\E )$, la stabilité de la surcohérence par foncteur de localisation
entraîne que
$\R \underline{\Gamma} ^\dag _{Z'} f ^! (\E (\hdag T))
\in F\text{-}D ^\mathrm{b} _\mathrm{surcoh} ( \D ^\dag _{\PP',\Q})$.
Ainsi, $\H ^j \R \underline{\Gamma} ^\dag _{Z'} f ^! (\E (\hdag T))$
est un $\D ^\dag _{\PP',\Q}$-module surcohérent.

Comme $\H ^j \E  |\U \riso \sp _{Y \hookrightarrow \U,+} (G_j)$
alors
$\R \underline{\Gamma} ^\dag _{Y'} g ^! (\H ^j \E  |\U ) \riso \sp _{Y' \hookrightarrow \U',+} (b^* G_j)$
(voir \cite[4.1.8]{caro-construction}).
Ainsi, pour tout entier $i \not =0$,
$\H ^i \R  \underline{\Gamma} ^\dag _{Y'} g ^! (\H ^j \E  |\U ) =0$.
Via la suite deuxième spectrale d'hypercohomologie de
$\R  \underline{\Gamma} ^\dag _{Y'} g ^! $ utilisée pour $\E |\U$, il en résulte
$\H ^j \R  \underline{\Gamma} ^\dag _{Y'} g ^! ( \E  |\U )
\riso
\R  \underline{\Gamma} ^\dag _{Y'} g ^! (\H ^j \E  |\U ) $.
Or, $\H ^j  \R \underline{\Gamma} ^\dag _{Z'} f ^! (\E (\hdag T)) |\U ' \riso
\H ^j  \R \underline{\Gamma} ^\dag _{Y'} g ^! (\E  |\U )$.
D'où :
$\H ^j  \R \underline{\Gamma} ^\dag _{Z'} f ^! (\E (\hdag T)) |\U ' \riso
\sp _{Y' \hookrightarrow \U',+} (b^* G_j)$.
Puisqu'en outre $\H ^j  \R \underline{\Gamma} ^\dag _{Z'} f ^! (\E (\hdag T))$ est
un $\D ^\dag _{\PP'} (\hdag T') _{\Q}$-module cohérent et à support dans le sous-schéma fermé lisse $Z'$,
il en découle,
avec la description de l'image essentielle du foncteur
$\sp _{Z'\hookrightarrow \PP', T',+}$ (voir \cite[4.1.9]{caro-construction} et \cite[2.2.12]{caro_courbe-nouveau},
ou la fin de la preuve de \ref{fibextfinidense2}),
qu'il existe un $F$-isocristal $E'_j$ surconvergent sur $Y'$
tel que
$ \H ^j \R \underline{\Gamma} ^\dag _{Z'} f ^! (\E (\hdag T))\riso \sp _{Z'\hookrightarrow \PP', T',+} (E'_j )$.
Comme
$\sp _{Z'\hookrightarrow \PP', T',+} (E'_j ) |\Y' \riso \sp _{Y' \hookrightarrow \U',+} (\widehat{E'} _j )$,
on en déduit alors
$\widehat{E'} _j \riso b ^* (G _j)$ (on rappelle que le foncteur
$\sp _{Y' \hookrightarrow \U',+} $ est pleinement fidèle).

\end{proof}

Or, pour $i\not =0$,
$\H ^i f_+ \sp _{Z'\hookrightarrow \PP', T',+} (E'_j )=0$.
En effet, comme $f_+ \sp _{Z'\hookrightarrow \PP', T',+} (E'_j )
\in D ^\mathrm{b} _\mathrm{coh} ( \D ^\dag _{\PP} (\hdag T) _{\Q})$
(car la cohérence est préservée par l'image directe d'un morphisme propre),
il suffit d'après \cite[4.3.12]{Be1} de le voir en dehors de $T$.
Or, comme $b$ est fini et étale, $f_+ \sp _{Z'\hookrightarrow \PP', T',+} (E'_j )|\U
\riso
g_+ \sp _{Y'\hookrightarrow \U', +} (\widehat{E} '_j )
\riso
\sp _{Y\hookrightarrow \U, +} (b _* \widehat{E} '_j )$,
où $\widehat{E} '_j$ est le $F$-isocristal convergent sur $Y'$ induit par $E' _j$
et
$b _* \smash{\widehat{E}} '_j$ est l'isocristal convergent sur $Y$ image directe par $b$ de $\widehat{E} '_j$
(ce dernier a bien un sens car $b$ est fini et étale, d'ailleurs en considérant les sections globales, il correspond canoniquement au foncteur oubli).
Via la deuxième suite spectrale d'hypercohomologie de $f _+$ utilisée pour le complexe
$\R \underline{\Gamma} ^\dag _{Z'} f ^! (\E (\hdag T))$,
d'après le lemme \ref{carac-surcoh-lemm*2},
il en résulte l'isomorphisme
$f _+  ( \H ^j \R \underline{\Gamma} ^\dag _{Z'} f ^! (\E (\hdag T)))\riso \H ^j f _+ (\R \underline{\Gamma} ^\dag _{Z'} f ^! (\E (\hdag T)))$.

On obtient alors par composition le morphisme :
$f _+  ( \H ^j \R \underline{\Gamma} ^\dag _{Z'} f ^! (\E (\hdag T)))
\rightarrow
\H ^j \E (\hdag T)$.
Pour terminer la {\og descente\fg},
nous construisons un second morphisme (qui donne une section de notre premier morphisme en dehors de $T$) de la façon suivante.
\medskip

Via l'équivalence de catégories (preuve de \cite[6.5.3]{caro_devissge_surcoh}) :
\begin{equation}
\label{653dev}
  F \text{-}\mathrm{Isoc} ^\dag( Y /K) \cong
F \text{-}\mathrm{Isoc} ^{\dag }(Y' /K) \times _{F \text{-}\mathrm{Isoc} ^{\dag } (Y ', Y' /K) }
F \text{-}\mathrm{Isoc} ^{\dag } (Y, Y /K),
\end{equation}
il existe un $F$-isocristal $E_j$ surconvergent sur $Y$ tel que $b^* (E _j) \riso E' _j$
et $\widehat{E} _j \riso G _j$.
D'après \cite[1.3.6]{caro-2006-surcoh-surcv}, on dispose de $\widetilde{\E}_j:= \sp _{Y, +}(E_j)\in
F \text{-}\mathrm{Isoc} ^{\dag \dag} ( Y /K) = F \text{-}\mathrm{Isoc} ^{\dag \dag} ( \PP, T, Z /K)$.
Lors de la preuve de \cite[6.3.1]{caro_devissge_surcoh},
nous avons construit (en {\og dualisant\fg} le morphisme d'adjonction)
le morphisme canonique
$\widetilde{\E} _j \rightarrow
f _+ \R \underline{\Gamma} ^\dag _{Z'} f ^! (\widetilde{\E} _j )$.
Or, $ \H ^j \R \underline{\Gamma} ^\dag _{Z'} f ^! (\E (\hdag T))\riso \sp _{Z'\hookrightarrow \PP', T',+} (E'_j )$
(voir le lemme \ref{carac-surcoh-lemm*2})
et
$\R \underline{\Gamma} ^\dag _{Z'} f ^! (\widetilde{\E} _j )\riso \sp _{Z'\hookrightarrow \PP', T',+} (b^* (E _j) )$
(d'après \cite[6.3.1]{caro_devissge_surcoh}).
Puisque $b^* (E _j) \riso E' _j$, on obtient alors
$\R \underline{\Gamma} ^\dag _{Z'} f ^! (\widetilde{\E} _j )\riso
 \H ^j \R \underline{\Gamma} ^\dag _{Z'} f ^! (\E (\hdag T))$.
 D'où la flèche :
$\widetilde{\E} _j \rightarrow
f _+   \H ^j \R \underline{\Gamma} ^\dag _{Z'} f ^! (\E (\hdag T))$.

En composant nos deux morphismes, on obtient :
$\widetilde{\E}_j \rightarrow \H ^j \E (\hdag T)$. Par construction, ce dernier est un isomorphisme en dehors
de $T$.
Comme c'est en outre un morphisme de
$\D ^\dag _{\PP} (\hdag T) _\Q$-modules cohérents, d'après \cite[4.3.12]{Be1}, on en déduit
$\widetilde{\E}_j \riso \H ^j \E (\hdag T)$.
En particulier,
$\H ^j \E (\hdag T)$ est un facteur direct de $f _+  \H ^j \R \underline{\Gamma} ^\dag _{Z'} f ^! (\E (\hdag T))$.

Comme $\R \underline{\Gamma} ^\dag _{Z'} f ^! (\E (\hdag T))\in F\text{-}D ^\mathrm{b} _\mathrm{surcoh} ( \D ^\dag _{\PP',\Q})$,
par stabilité de la surcohérence par l'image directe d'un morphisme propre on obtient alors
$f _+   \H ^j \R \underline{\Gamma} ^\dag _{Z'} f ^! (\E (\hdag T))
\in F\text{-}D ^\mathrm{b} _\mathrm{surcoh} ( \D ^\dag _{\PP,\Q})$.
Il en résulte que $\H ^j \E (\hdag T)\in F\text{-}D ^\mathrm{b} _\mathrm{surcoh} ( \D ^\dag _{\PP,\Q})$, ce qui est équivalent à
dire que $\E (\hdag T)\in F\text{-}D ^\mathrm{b} _\mathrm{surcoh} ( \D ^\dag _{\PP,\Q})$.

Via le triangle de localisation de $\E$ en $T$ (voir \ref{tri-loc-gen}),
pour vérifier que $\E \in F\text{-}D ^\mathrm{b} _\mathrm{surcoh} ( \D ^\dag _{\PP,\Q})$,
il suffit d'établir que
$\R \underline{\Gamma} ^\dag _{T} (\E ) \in F\text{-}D ^\mathrm{b} _\mathrm{surcoh} ( \D ^\dag _{\PP,\Q})$.
Si $\dim  T \cap X \leq n-1$ (ce qui est le cas par exemple si $Z =X$), cela résulte de l'hypothèse que $(Q _{n-1})$ est vraie.
Autrement, si $\dim  T \cap X = \dim X =n$, on remarque que $T \cap X$ possède au moins une composante
irréductible de dimension égale à $n$ de moins que $X$. On conclut alors par récurrence sur
le nombre de composantes irréductibles de dimension égale à $n$.
\end{proof}

\begin{rema}
  L'hypothèse que $\PP$ soit propre et la structure de Frobenius dans \ref{carac-surcoh} sont sans-doute superflues.
  Elles nous ont cependant été utiles dans la preuve via le corollaire \ref{fibextfinidense}, l'équivalence de catégories
  \ref{653dev} et aussi via l'utilisation de l'équivalence de catégories induite par $\sp _{Y+}$ (voir \ref{eqcat-gen}).
Une fois étendues ces trois résultats au cas des isocristaux partiellement surconvergents (ce qu'il est très loin d'être évident),
on obtiendrait alors l'extension de \ref{carac-surcoh} du cas où $\PP$ est propre et lisse avec structure de Frobenius
au cas où $\PP$ est séparé et lisse sans structure de Frobenius.
\end{rema}

\bibliographystyle{smfalpha}
\bibliography{bib1}

\bigskip
\noindent Daniel Caro\\
Arithmétique et Géométrie algébrique, Bât. 425\\
Université Paris-Sud\\
91405 Orsay Cedex\\
France.\\
email: daniel.caro@math.u-psud.fr

\end{document}